\documentclass[12pt]{article}
\usepackage{amssymb}
\usepackage{amsmath}
\usepackage{graphicx}
\usepackage{epstopdf}
\usepackage{subfigure}
\usepackage{longtable}
\usepackage{url}
\usepackage[hyphenbreaks]{breakurl}
\usepackage[colorlinks,linkcolor=red]{hyperref}

\numberwithin{equation}{section}

\def\kasten{$~~\mbox{\hfil\vrule height6pt width5pt depth-1pt}$ }

\newtheorem{Theorem}{Theorem}[section]

\newtheorem{Proposition}[Theorem]{Proposition}

\newtheorem{Example}[Theorem]{Example}
\newtheorem{Remark}[Theorem]{Remark}

\begin{document}

\def\kasten{$~~\mbox{\hfil\vrule height6pt width5pt depth-1pt}$ }
\newtheorem{theorem}{Theorem}

\centerline{\large \bf Discovering transition phenomena from data of }
\centerline{\large \bf stochastic dynamical systems  with L\'evy noise }

\bigskip
\centerline{\bf Yubin Lu$^{a,}\footnote{yubin\_lu@hust.edu.cn}$ and
Jinqiao
Duan$^{b,}\footnote{Corresponding author: duan@iit.edu}$}
\smallskip
\centerline{${}^a$ School of Mathematics and Statistics \& Center for Mathematical Sciences}
\smallskip
\centerline{Huazhong University of Science and Technology, Wuhan 430074, China}
\smallskip
\centerline{${}^b£¬$ Department of Applied Mathematics, College of Computing}
\smallskip
\centerline{Illinois Institute of Technology, Chicago, IL 60616, USA}
\centerline{August 14, 2020}

\begin{abstract}
It is a challenging issue to analyze complex dynamics  from observed and simulated data.  An advantage of extracting dynamic behaviors from data is that this approach enables the investigation of nonlinear phenomena whose mathematical models are unavailable. The purpose of this present work  is to extract information about   transition phenomena (e.g.,  mean exit time and escape probability), from data of  stochastic differential  equations with non-Gaussian L\'evy noise. As a tool in  describing dynamical systems, the Koopman semigroup  transforms a nonlinear system  into a linear system, but at the cost of elevating a finite dimensional problem into an infinite dimensional one. In spite of this, using the relation between the stochastic Koopman semigroup and the infinitesimal generator of a stochastic differential equation, we learn the  mean exit time and escape probability  from data.  Specifically, we  first obtain  a finite dimensional  approximation of the infinitesimal generator     by   an  extended dynamic mode decomposition algorithm.   Then we     identify  the  drift coefficient, diffusion coefficient and anomalous diffusion coefficient for the stochastic differential equation. Finally,  we compute the mean exit time and escape probability by finite difference discretization of  the associated  nonlocal partial differential equations. This approach  is applicable in extracting transition information from data of  stochastic differential equations   with either (Gaussian) Brownian motion  or (non-Gaussian)  L\'evy motion.  We  present  one - and two-dimensional examples to demonstrate the effectiveness of our approach.

\end{abstract}

\bigskip\par

\textsl{The transition behavior in stochastic dynamical systems is often quantified  in terms of deterministic indexes (e.g., mean exit time, escape probability), leading to simple and effective prediction. However, the discovery of such indexes usually relies on a deep understanding and priori knowledge of the system, as well as extensive and time-consuming mathematical justification. Recent developments in data-driven computational approaches suggest alternative ways toward discovering deterministic indexes directly from data. It is intriguing to note that  obtaining such deterministic indexes for stochastic dynamical systems accurately and efficiently is still a significant challenge. In this paper, we extract deterministic indexes from data of  a stochastic dynamical system with L\'evy noise,  after identifying the  system by  adapting a method to learn  the associated  Koopman semigroup and its generator.  }

\section{Introduction}

Stochastic differential equations  are widely used to model  complex phenomena in climate dynamics, molecular dynamics,   biophysical dynamics,  and  other engineering systems, under random fluctuations. Usually, we build  a mathematical model for such a    system based on   available  governing laws,  and then analyze or simulate the model to gain insights about its nonlinear phenomena (e.g., transition phenomena \cite{Duan}). However,  sometimes it is difficult to  build a model for a  complex  phenomenon due to lack of scientific understanding.  At some other times, a mathematical model based on governing laws may be too complicated for analysis.  Fortunately, with the development of observing  techniques  and   computing power, there are  a lot of valuable observation data when a model is unavailable,  or  simulated data when a model is too complicated yet computable. Therefore, it is desirable  to discover  dynamical quantities and  gain insights about stochastic phenomena directly from these data sets.  \\

Data-driven analysis of complex dynamic behavior has recently attracted increasing attention.   Many authors have proposed  various methods to analyze complex dynamic behaviors from data. These include, for example, stochastic parametrization \cite{AJLu}, parameters estimation of stochastic differential equations using Gaussian process \cite{GPs1,GPs2}, Sparse Identification of Nonlinear Dynamics \cite{SINDy1,SINDy2}, Kramers-Moyal expansion \cite{KM}, the Koopman operator analysis,  for analyzing complex dynamic behaviors \cite{DMD,EDMD,SKO1,SKO2,KO,SKO3}, and transfer operators for analyzing stochastic dynamic behaviors \cite{Dellnitz,Froyland,Klus1,Klus2,Metzner,Tantet,Thiede}, among others. However, most of these authors focused on deterministic dynamical systems, or stochastic dynamical systems with Gaussian noise.

In this present work, we devise a method to extract the mean exit time and escape probability from data, for a stochastic dynamical system with non-Gaussian stable L\'evy noise.  These two quantities provide insights about transition phenomena under the interaction of nonlinearity and uncertainty \cite{Duan}.  We first infer the generator for this non-Gaussian stochastic dynamical system, with help of   the Koopman operator, which fortunately also deals with      dynamical  systems with noise.  Then we identify the coefficients in this stochastic system, and finally we discretize nonlocal partial differential equations to obtain the mean exit time and escape probability.  The main contribution of this present work is an algorithm to estimate the coefficient of the non-Gaussian Levy noise from dynamical data, and thus further enable to extract transition information from the same data. \\

The remainder of this paper is arranged  as follows: In Section 2, we present the theoretical background, including the definition of Koopman operator and   infinitesimal generator, and then the mean exit time and escape probability. In Section 3, we will introduce the adapted  EDMD algorithm and the method of system identification. In Section 4, we will  conduct numerical experiments,  including stochastic differential equations with  Brownian motion and L\'evy motion, in   one-dimensional and two-dimensional cases.  We leave the discussion    in the final  section. \\

\section{Theory}

We recall the Koopman semigroup of operators and its generator in Section \ref{Koopman}, and then introduce  deterministic quantities, mean exit time and escape probability,  for stochastic dynamics in Section \ref{quantities}.

\subsection{The Koopman operator and the  generator}\label{Koopman}

Let $\mathcal{M}$ be the state space in $\mathbb{R}^d$. Consider a stochastic differential equation in $\mathcal{M}$
\begin{align}
\label{SDE1}
dX_{t}=b(X_{t})dt+\sigma_{1}(X_{t})dW_{t}+\sigma_{2}d\tilde{L}_{t}^\alpha,
\end{align}
with  the drift coefficient (or vector field) $b$ in $\mathbb{R}^d$, diffusion coefficient  $d \times d$ diagonal matrix $\sigma_1$,   $d-$dimensional Brownian motion $W_t$ and   $d-$dimensional L\'evy motion $\tilde{L}_{t}^\alpha$.   The components of $\tilde{L}_{t}^\alpha$ are one-dimensional mutually independent and symmetric with the same triple $(0,0,\tilde{\nu})$ L\'evy motion, where $\tilde{\nu}(dx)=C_{\alpha}I_{\{\Vert x\Vert<c\}}(x)\frac{dx}{|x|^{1+\alpha}}$ for $x\in \mathbb{R}\setminus\{0\}$. Here $c$ is a positive constant (usually we take $c=1$), $I$ is a indicator function and
$$C_{\alpha}=\frac{\alpha\Gamma(\frac{1+\alpha}{2})}{2^{1-\alpha}\pi^{\frac{1}{2}}\Gamma(1-\frac{\alpha}{2})}.$$
 Moreover,  the anomalous diffusion coefficient  $\sigma_2$ is a constant $d \times d$ diagonal matrix. Since L\'evy motion is symmetric, without loss of generality, we can set $\sigma_{2}\geq0$. In this paper we assume that the stability index $\alpha$ is known. More details about L\'evy motion are in Appendix. \\
The associated  Koopman semigroup of operators $\{\mathcal{K}^t\}$ is defined as
\begin{align}
(\mathcal{K}^tf)(x)=\mathbb{E}[f(\Phi^t(x,\omega))],
\end{align}
where $f:\mathbb{R}^d \rightarrow \mathbb{R}$ a real-valued measurable function, $\Phi^t(x,\omega)$ is the solution of (\ref{SDE1}) with initial value $x$ and $\mathbb{E}[\cdot]$ denotes the expected value.  The infinitesimal generator $\mathcal{L}$ of this Koopman semigroup is the derivative of $\mathcal{K}^t$ at $t=0$.
\begin{align}
\mathcal{L}f=\lim\limits_{t\to0}\frac{1}{t}(\mathcal{K}^tf-f).
\end{align}
Thus, if $f\in Dom(\mathcal{L})$, the infinitesimal generator  is a nonlocal partial differential operator
\begin{align}
\label{generator1}
\mathcal{L}f&=\sum_{i=1}^d b_{i}\frac{\partial{f}}{\partial{x_{i}}}+\frac{1}{2}\sum_{i=1}^d\sum_{j=1}^d a_{ij}\frac{\partial^{2}f}{\partial{x_{i}}\partial{x_{j}}} \nonumber\\
&+C_{\alpha}\sum_{i=1}^d\int_{\mathbb{R}\backslash\{0\}} [f(x_{i}+(\sigma_{2})_{ii}y_{i})-f(x_{i})]I_{\{\Vert y_{i}\Vert<c\}}(y_{i})\, \frac{dy_{i}}{|y_{i}|^{1+\alpha}},
\end{align}
where  the   coefficient $a=\sigma_{1}\mathbf{\sigma}_{1}^\top$ is sometimes also called   diffusion coefficient.   For more details see \cite{Chen,Sun}.  \\

When the L\'evy motion is absent, the third term in this  generator is absent (because Brownian motion has no jumps) and the generator becomes a local partial differential operator
 \begin{align}
\label{generator2}
\mathcal{L_B}f&=\sum_{i=1}^d b_{i}\frac{\partial{f}}{\partial{x_{i}}}+\frac{1}{2}\sum_{i=1}^d\sum_{j=1}^d a_{ij}\frac{\partial^{2}f}{\partial{x_{i}}\partial{x_{j}}}.
\end{align}

\subsection{Deterministic quantities for  transition phenomena }\label{quantities}
We mainly discuss mean exit time and escape probability with non-Gaussian noise,  as deterministic quantities that carry dynamical information for solution orbits of (\ref{SDE1}). In what follows, we assume that $D$  is a regular bounded open domain in $\mathbb{R}^d$. \\
\\
\textbf{Mean Exit Time} \\
\\
We first consider the mean exit time. Define the first exit time $\tau_{x}$ of a solution path starting at $x$ from domain $D$ as
\begin{align}
\tau_{x}(\omega)=inf\{t\geq 0, X_{0}=x, X_{t}(\omega,x) \not\in D\}.
\end{align}
The mean exit time is then denoted by
\begin{align}
u(x)=\mathbb{E}^x[\tau_{x}(\omega)],
\end{align}
for $x \in D$. \\
As usual, we assume that the drift coefficient and diffusion coefficient satisfy an appropriate local Lipschitz condition. Moreover, we assume that the generator operator of this system is uniformly elliptic. For more details, see \cite{Duan,Applebaum}. Under these conditions, the mean exit time $u(x)$, for a solution path starting at $x \in D$, satisfies the following nonlocal partial differential equation\\
\begin{align}\label{eqn:MET1}
\mathcal{L}u=-1,
\end{align}
with ``an external'' Dirichlet boundary condition
\begin{align}\label{eqn:MET2}
u\mid_{D^c}=0,
\end{align}
where $\mathcal{L}$ is the generator (\ref{generator1}), and $D^c$ is the complement of the domain $D$ in $\mathbb{R}^d$. \\
\\
\textbf{Escape Probability} \\
\\
Here we just recall  the escape probability  for stochastic system with  L\'evy motion.
Assume that $X_{t}(\omega)$ is the solution of (\ref{SDE1}). Now we consider the escape probability of $X_{t}(\omega)$. Define the first exit time
\begin{align}
\tau_{D^c}(\omega)\triangleq  \inf\{t>0: X_{t}(\omega) \in D^c\}.
\end{align}
We take a subset $U$ of $D^c$, and define the likelihood that $X_{t}(\omega)$ exits firstly from $D$ by landing in the target set $U$ as the escape probability from $D$ to $U$, denoted by $p(x)$. That is
\begin{align}
p(x)=\mathbb{P}\{X_{\tau_{D^c}}(\omega) \in U\}.
\end{align}
Then the escape probability $p$, for the dynamical system driven by L\'evy motion (\ref{SDE1}), from $D$ to $U$, is the solution of the following nonlocal partial differential equations with Balayage-Dirichlet external boundary conditions
\begin{align}\label{eqn:EP}
             &\mathcal{L}p=0, \nonumber\\
             &p\mid_{U}=1, \nonumber\\
             &p\mid_{D^c\setminus U}=0,
\end{align}
where $\mathcal{L}$ is the generator (\ref{generator1}). \\

For a system  with only Brownian motion,  the sample  paths are continuous and thus can only exit  domain $D$ by passing  its  boundary $\partial D$, or escape through a portion $U$ of the boundary. Therefore, the mean exit time    $u$ satisfies
  \begin{align}\label{eqn:MET3}
\mathcal{L_B} \; u=-1, \;\;\;   u\mid_{\partial D}=0.
\end{align}
Moreover, the escape probability $p$, through a portion  $U$ of the boundary, satisfies
\begin{align}\label{eqn:EP3}
             &\mathcal{L_B}\; p=0, \nonumber\\
             &p\mid_{U}=1, \nonumber\\
             &p\mid_{\partial D \setminus U}=0.
\end{align}
  For more details see \cite{Duan}. \\

\section{Numerical schemes}
In this section, we will review briefly the extended dynamic mode decomposition in Section \ref{secEDMD},  and then identify system in Section  \ref{SI}. Finally, we will propose a algorithm for discovering transition phenomena from data in Section \ref{algorithm}.

\subsection{Extended Dynamic Mode Decomposition}\label{secEDMD}
In this subsection, we review breifly extended dynamic mode decomposition(EDMD), which can approximate Koopman operators. Of course, EDMD can also approximate Koopman eigenvalue, eigenfunction and mode tuples as described by \cite{EDMD}, but we will not consider them here.\\

We adopt the notations from \cite{EDMD}. There are two prerequisites in the EDMD procedure: \\
(a) A data set of snapshot pairs, i.e., $\{(x_{m}, y_{m})\}_{m=1}^M$ that we will organize as a pair of data sets,
\begin{align}
X=[x_{1}, x_{2},\ldots,x_{M}], \nonumber\\
Y=[y_{1}, y_{2},\ldots,y_{M}],
\end{align}
where $x_{i}\in\mathcal{M}$ and $y_{i}\in\mathcal{M}$ satisfy $y_{i}=\Phi^t(x_{i},\omega)$.\\
(b) Assume that $\mathcal{F}=L^2(\mathcal{M})$. A dictionary of observables $\mathcal{D}=\{\psi_{1}, \psi_{2}, \ldots, \psi_{N_{k}}\}$, where $\psi_{i}\in\mathcal{F}$, whose span we denote as $\mathcal{F}_{\mathcal{D}}\subset\mathcal{F}$; for the sake of simplicity, we also denote it by a vector-valued function $\Psi:\mathcal{M} \rightarrow \mathbb{C}^{1\times N_{K}}$ as follows
\begin{align}
\Psi(x)=[\psi_{1}(x), \psi_{2}(x), \ldots, \psi_{N_{k}}(x)].
\end{align}
A function $f\in \mathcal{F}_{\mathcal{D}}$ can be written as $f = B^\top\Psi^\top$, where $B$ is a weight vector. Because $\mathcal{F}_{\mathcal{D}}$ is not an invariant subspace of $\mathcal{K}$, we obtain the follow relation with a residual term $r(x)$:
\begin{align}
(\mathcal{K}^{\Delta t}f)(x)&=\mathbb{E}[f\circ\Phi^{\Delta t}(x,\omega)] \nonumber\\
&=B^\top\mathbb{E}[\Psi^{\top}\circ\Phi^{\Delta t}(x,\omega)]+r(x) \nonumber\\
&=(KB)^\top\Psi^{\top}(x)+r(x).
\end{align}
Then we will minimize  the residual term $r(x)$ to obtain the matrix $K=G^+A$, where $+$ denotes the pseudoinverse and
\begin{align}
G=\frac{1}{M}\sum_{m=1}^{M}\Psi(x_{m})^*\Psi(x_{m}),
\end{align}
\begin{align}
A=\frac{1}{M}\sum_{m=1}^{M}\Psi(x_{m})^*\Psi(y_{m}).
\end{align}
Here  matrices  $K,G,A\in\mathbb{C}^{N_{K}\times N_{K}}$. \\
Therefore, we have  a matrix $K$ such that
\begin{align}
(\mathcal{K}^{\Delta t}f)(x)\approx (KB)^\top\Psi^\top,
\end{align}
where $\Delta t$ is small enough.\\
Using the definition for  the generator $\mathcal{L}f=\lim\limits_{t\to0}\frac{1}{t}(\mathcal{K}^tf-f)$, we can get a finite-dimensional approximation of generator $\mathcal{L}$ by finite-dimensional approximation $K$ of the stochastic Koopman operator $\mathcal{K}$. This approximated generator is denoted by $L$. \\

This approximated generator $L$  can not  be directly used     in solving     partial differential equations for mean exit time and escape probability, as these two dynamical quantities    satisfy special boundary conditions as in Section \ref{quantities}.  So we next need to identify coefficients in the stochastic differential  equation (\ref{SDE1}), and use these coefficients in discretizing partial differential equations to obtain  mean exit time and escape probability.

\subsection{System identification}\label{SI}

In this subsection, we will state the method for identifying system parameters using generator operator. Assume that the drift coefficient and diffusion coefficients of a stochastic differential equation can be expanded by basis functions $\Psi(x)$.  We assume that the diffusion coefficients $\sigma_{1}$ and $\sigma_{2}$ are diagonal matrices. Indeed, $\sigma_{1}$ does not need to be a diagonal matrix, but our method requires that $\sigma_{2}$ to be diagonal; See Remark \ref{remark}.\\
Given the full-state observable
$$f_{1i}(x)=x_{i}, i=1,2,\ldots,d,$$
where $x=(x_1,x_2,\ldots,x_d)\in \mathbb{R}^d$. \\
With the aid of generator operator $\mathcal{L}$ and its finite-dimensional approximation $L$, described in previous section, it is possible to discover the governing equations. Assume that $B_{1i}\in\mathbb{R}^{N_{K}\times 1}$ is the vector such that $f_{1i}(x)=B_{1i}^\top\Psi(x)^\top$. We can express the drift coefficient in terms of the basis functions, i.e.,
\begin{align}\label{eqn:SI1}
(\mathcal{L}f_{1})(x)=b(x)\approx (LB_{1i})^\top\Psi(x)^\top.
\end{align}
Similarity, let $f_{2i}(x)=x_{i}^2$, we can get the diffusion coefficients by finding another matrix $B_{2i}\in\mathbb{R}^{N_{K}\times 1}$ such that $f_{2i}(x)=B_{2i}^\top\Psi(x)^\top$
\begin{align}\label{eqn:SI2}
(\mathcal{L}f_{2i})(x)&=2b_{i}(x)x_{i}+a_{ii}(x) \nonumber\\
&+C_{\alpha}\int_{\mathbb{R}\backslash\{0\}} \frac{[(x_{i}+(\sigma_{2})_{ii}y_{i})^2-x_{i}^2]I_{\{\Vert y_{i}\Vert<c\}}(y_{i})}{|y_{i}|^{1+\alpha}}\, dy_{i} \approx (LB_{2i})^\top\Psi(x)^\top,
\end{align}
Note that, in order to make sure the integral in (\ref{eqn:SI2}) exists, we assume that the L\'evy motion has bounded jump size with bound $c$ (a positive constant). Moreover, using the properties of odd functions, the integral term can be represented  by \\
\begin{align}\label{eqn:SI3}
\int_{\mathbb{R}\backslash\{0\}} \frac{[(x_{i}+(\sigma_{2})_{ii}y_{i})^2-x_{i}^2]I_{\{\Vert y_{i}\Vert<c\}}(y_{i})}{|y_{i}|^{1+\alpha}}\, dy_{i}=(\sigma_{2})_{ii}^2\int_{\mathbb{R}\backslash\{0\}} \frac{y_{i}^2I_{\{\Vert y_{i}\Vert<c\}}(y_{i})}{|y_{i}|^{1+\alpha}}\, dy_{i}.
\end{align}
Therefore we can rewrite (\ref{eqn:SI2}) as follows
\begin{align}\label{eqn:SI4}
(\mathcal{L}f_{2i})(x)&=2b_{i}(x)x_{i}+a_{ii}(x) \nonumber\\
&+(\sigma_{2})_{ii}^{2}C_{\alpha}\int_{\mathbb{R}\backslash\{0\}} \frac{y_{i}^{2}}{|y_{i}|^{1+\alpha}}I_{\{\Vert y_{i}\Vert<c\}} (y_{i})\, dy_{i}.
\end{align}
Note that there are two unknown coefficients in (\ref{eqn:SI4}), i.e., diffusion coefficient $a$ and anomalous diffusion coefficient $\sigma_{2}$. In general,  one equation can not determine two unknown coefficients. However, we will explain  this can actually be done in our particular setting.  We present this for a scalar system     in    the  following  proposition. \\

\begin{Proposition}\label{Prop}
Consider a scalar stochastic differential equation
\begin{align}\label{SDE4}
dX_{t}=b(X_{t})dt+\sigma_{1}(X_{t})dW_{t}+\sigma_{2}d\tilde{L}_{t}^\alpha,
\end{align}
where the drift coefficient $b(x)=\sum_{i=0}^{p_{1}}\xi_{i}x^i$,  diffusion coefficient $\sigma_{1}(x)=\sum_{j=0}^{p_{2}}\eta_{j}x^j$,  and  anomalous  diffusion coefficient $\sigma_{2}$ is a constant. Let $\Psi(x)=[1, x, \ldots, x^{N_{k}}]$ be the polynomial basis, with $N_{k}$   large enough. If $\sigma_{1}$ is not a constant and the equation (\ref{eqn:SI6}) is well defined, then $\sigma_{1}$ and $\sigma_{2}$ can be separated by (\ref{eqn:SI4}).\\

\textbf{Proof} \quad Using the polynomial basis $\Psi(x)$, we   express  the coefficients $b(x)=\xi^\top\Psi(x)^\top$ and $\sigma_{1}(x)=\eta^\top\Psi(x)^\top$, with coefficients  $\xi=[\xi_{0},\xi_{1},\ldots,\xi_{p_{1}},0,\ldots,0]^\top$ and $\eta=[\eta_{0},\eta_{1},\ldots,\eta_{p_{2}},0,\ldots,0]^\top$. Note that we   can also  express  $\sigma_{2}$ using the polynomial basis, i.e., $\sigma_{2}=\zeta^\top\Psi(x)^\top$, with coefficient $\zeta=[\sigma_{2},0,\ldots,0]^\top$.\\
From (\ref{eqn:SI1}),   we have the coefficient   $\xi=LB_{1}$, i.e., $\xi$ is known. Furthermore, note that
$\sigma_{1}^{2}(x)=\tilde{\eta}^\top\Psi(x)^\top$ and $\sigma_{2}^{2}=\tilde{\zeta}^\top\Psi(x)^\top$,
where
$$\tilde{\eta}=[\tilde{\eta}_{0}, \tilde{\eta}_{1}, \tilde{\eta}_{2}, \ldots, \tilde{\eta}_{2p_{2}}, 0, \ldots, 0]^\top$$
and
$$\tilde{\zeta}=[\sigma_{2}^2,0,\ldots,0]^\top.$$
It should be noted that $\tilde{\eta}_{j}$ is the coefficient of $x^j$ for $\sigma_{1}^{2}(x)$, $j=0, 1, 2,\ldots, 2p_{2}$.\\
Therefore, from (\ref{eqn:SI4}) we have
\begin{align}
2(\xi)_{+}+\tilde{\eta}+\tilde{C}\tilde{\zeta}=LB_{2},
\end{align}
where $\tilde{C}=C_{\alpha}\int_{\mathbb{R}\backslash\{0\}} \frac{y^2}{|y|^{1+\alpha}}I_{\{\Vert y\Vert<c\}} (y)\, dy$ is a constant and the symbol $(v)_{+}$ for column vector $v=[v_{1},\ldots,v_{n}]^\top$ is a shift operator. i.e., $(v)_{+}=[0,v_{1},\ldots,v_{n-1}]^\top$.
Denoting  $\rho\triangleq LB_{2}-2(\xi)_{+}$, we have
\begin{align}
\label{eqn:SI5}
\tilde{\eta}+\tilde{C}\tilde{\zeta}=\rho,
\end{align}
or
\begin{align}
\label{eqn:SI6}
 \begin{pmatrix} 
     \tilde{\eta}_{0}  \\
     \tilde{\eta}_{1}  \\
     \vdots \\
     \tilde{\eta}_{2p_{2}}  \\
     0 \\
     \vdots \\
     0
 \end{pmatrix} 
 +
  \begin{pmatrix} 
     \tilde{C}\sigma_{2}^2  \\
     0  \\
     \vdots \\
     0  \\
     0 \\
     \vdots \\
     0
 \end{pmatrix} 
 =
  \begin{pmatrix} 
     \rho_{0}  \\
     \rho_{1}  \\
     \vdots \\
     \rho_{2p_{2}}  \\
     0 \\
     \vdots \\
     0
 \end{pmatrix}.
\end{align}
Our goal is to separate $\sigma_{1}^2$ and $\sigma_{2}^2$, i.e., we only need to determine the $2p_{2}+2$ unknown coefficients $\tilde{\eta}_{0}, \tilde{\eta}_{1}, \ldots, \tilde{\eta}_{2p_{2}}$ and $\sigma_{2}^2$. Since       the equation (\ref{eqn:SI6}) is well defined and           $\tilde{\eta}$ consists of coefficients $\eta_{0}, \eta_{1}, \ldots, \eta_{p_{2}}$, we can separate $\sigma_{1}^2$ and $\sigma_{2}^2$ by a backward iteration. More precisely, first, from (\ref{eqn:SI6}), we know that $\tilde{\eta}_{i}= \rho_{i}, i=1,2,\ldots,2p_{2}$. Second, as for $\tilde{\eta}_{0}$, note that  $\tilde{\eta}_{0}=\eta_{0}^2$ and realize that $\{\tilde{\eta}_{i}\}_{i=1}^{2p_{2}}$ are expressed in terms of  variables  $\eta_{0}, \eta_{1},\ldots,\eta_{p_{2}}$. The equation (\ref{eqn:SI6}) is well defined, which ensures the uniquely solvability for equations ${\tilde{\eta}_{i}}=\rho_{i}, i=1,2,\ldots,2p_{2}$. Furthermore, $\sigma_{1}$ is not a constant, i.e., $p_{2}\geq 1$, which ensures the number of these equations is large than the number of variables (recalling that we have $2p_{2}$ equations, $p_{2}+1$ variables and the relationship $2p_{2}\geq p_{2}+1$). Therefore, we can obtain $\eta_{0}^2$ from such $2p_{2}$ equations. It means that the coefficient $\tilde{\eta}_{0}$ is known now and this determines  $\sigma_{1}^2$.  Finally, by the first component of  equation (\ref{eqn:SI6}), we obtain $\sigma_{2}^2=\frac{\rho_{0}-\tilde{\eta}_{0}}{\tilde{C}}$.   That is why we call it as a backward iteration. This completes the proof.
$\hfill\blacksquare$

\end{Proposition}

\begin{Remark}\label{remark2}
In this scalar setting of  Proposition \ref{Prop}, we can actually only determine $\sigma_{1}^2$ via its coefficients $\tilde{\eta}$.  Also, after determining $\sigma_{2}^2$ we uniquely get $\sigma_2$, which is assumed to be non-negative.   Note that, as seen in the generator (\ref{generator1}),  $\sigma_{1}^2$  and  $\sigma_2$ are  just what we need.

\end{Remark}

\begin{Remark}
In Proposition \ref{Prop},  the equation  (\ref{eqn:SI6}) is a reformulation of (\ref{eqn:SI4}) in the sense of approximation. The equation (\ref{eqn:SI6}) should be well defined when $L$ is a good approximation of $\mathcal{L}$.
\end{Remark}

Here we present a simple example to show how to separate coefficients $\sigma_{1}^2$  and  $\sigma_2^2$. \\
\begin{Example}
Consider a scalar stochastic differential equation
\begin{align}
dX_{t}=b(X_{t})dt+\sigma_{1}(X_{t})dW_{t}+\sigma_{2}d\tilde{L}_{t}^\alpha,
\end{align}
where $b(x)=\sum_{i=0}^m \xi_{i}x^i$, $\sigma_{1}(x)=\eta_{1}x+\eta_{0}$ and $\sigma_{2}$ is a constant. Let $\Psi(x)=[1, x, \ldots, x^{N_{k}}]$, where $N_{k}$ is large enough. Therefore, we can rewrite the coefficients $b(x)=\xi^\top\Psi(x)^\top$ and $\sigma_{1}(x)=\eta^\top\Psi(x)^\top$, where $\xi=[\xi_{0},\xi_{1},\ldots,\xi_{m},0,\ldots,0]^\top$ and $\eta=[\eta_{0},\eta_{1},0,\ldots,0]^\top$. Note that we also can rewrite $\sigma_{2}$ using polynomial basis, i.e., $\sigma_{2}=\zeta^\top\Psi(x)^\top$, where $\zeta=[\sigma_{2},0,\ldots,0]^\top$.\\
From (\ref{eqn:SI1}), we know that $\xi=LB_{1}$. Furthermore, note that $\sigma_{1}^{2}(x)=\tilde{\eta}^\top\Psi(x)^\top$ and $\sigma_{2}^{2}=\tilde{\zeta}^\top\Psi(x)^\top$, where $\tilde{\eta}=[\eta_{0}^2,2\eta_{0}\eta_{1},\eta_{1}^2,0,\ldots,0]^\top$ and $\tilde{\zeta}=[\sigma_{2}^2,0,\ldots,0]^\top$. Therefore, from (\ref{eqn:SI4}) we have

\begin{align}
2(\xi)_{+}+\tilde{\eta}+\tilde{C}\tilde{\zeta}=LB_{2}.
\end{align}
Let $\rho\triangleq LB_{2}-2(\xi)_{+}$, we have

\begin{align}
\label{eqn:SI5}
\tilde{\eta}+\tilde{C}\tilde{\zeta}=\rho.
\end{align}

or
\begin{align}
\label{eqn:SI7}
 \begin{pmatrix} 
     \eta_{0}^2  \\
     2\eta_{0}\eta_{1}  \\
     \eta_{1}^2 \\
     0 \\
     \vdots \\
     0
 \end{pmatrix} 
 +
  \begin{pmatrix} 
     \tilde{C}\sigma_{2}^2  \\
     0  \\
     0  \\
     0 \\
     \vdots \\
     0
 \end{pmatrix} 
 =
  \begin{pmatrix} 
     \rho_{0}  \\
     \rho_{1}  \\
     \rho_{2}  \\
     0 \\
     \vdots \\
     0
 \end{pmatrix}
\end{align}
From (\ref{eqn:SI7}), we have $\eta_{1}^2=\rho_{2}$, $\eta_{0}^2=(\frac{\rho_{1}}{2\eta_{1}})^2$ and $\sigma_{2}^2=\frac{\rho_{0}-\eta_{0}^2}{\tilde{C}}$. In other words, the coefficients $\sigma_{1}^2$ and $\sigma_{2}^2$ are separated. Here we assume that the stable index $\alpha$ is known.

\end{Example}

\begin{Remark}\label{remark}
If $\sigma_{1}$ is not a diagonal matrix, we can take $f(x)=x_{i}x_{j}$ to obtain $a_{ij}$, $x\in \mathbb{R}^d$.
\end{Remark}

\subsection{Algorithm}\label{algorithm}
We implement the algorithm for discovering dynamics from data in the MATLAB programming language. The basic algorithm is outlined in Algorithm 1. The data that support the findings of this study are openly available in GitHub, reference number \cite{Code}.
\begin{table}[h]
\label{table:Algorithm1}
        \begin{tabular}{l}
        \hline
        \textbf{Algorithm 1} Discovering dynamics from data\\
        \hline
        Require: Basis functions $\Psi(x)=[\psi_{1}(x), \psi_{2}(x), \ldots, \psi_{N_{k}}(x)]$.\\
        INPUT: A data set of snapshot pairs $\{(x_{m}, y_{m})\}_{m=1}^M$.\\
        \\
        1. Setting $G=\frac{1}{M}\sum_{m=1}^{M}\Psi(x_{m})^*\Psi(x_{m})$ and $A=\frac{1}{M}\sum_{m=1}^{M}\Psi(x_{m})^*\Psi(y_{m})$.\\
        2. Setting $K=G^+A$.\\
        3. System identification by (\ref{eqn:SI1}) and (\ref{eqn:SI4}).\\
        4. Using a finite difference method to obtabin mean exit time and \\
        \quad escape probability by solving  (\ref{eqn:MET1}), (\ref{eqn:MET2}) and (\ref{eqn:EP}).\\
        \\
        OUTPUT: Mean exit time and escape probability.\\
        \hline
       \end{tabular}
\end{table}

\begin{Remark}
In order to simulate the numerical solution of  nonlocal partial differential equations  (\ref{eqn:MET1}), (\ref{eqn:MET2}) and (\ref{eqn:EP}), we use a finite difference method from our earlier works \cite{Gao,Chen}. \\
\end{Remark}

\section{Numerical experiments }
In this section we will present  a few  examples to illustrate  our method for extracting mean exit time and escape probability, from simulated data sets of stochastic dynamical systems.  We first consider    stochastic differential equations  driven by Brownian motion, and then we consider   stochastic differential equations     driven by L\'evy motion. Moreover,  we will illustrate that the method is valid for one-dimensional or two-dimensional cases. In the following, we choose the polynomials as basis functions.\\

\subsection{Dynamical systems with Brownian  motion  }

\begin{Example} Consider a scalar double-well system
\begin{align}
\label{SDE2}
dX_{t}=(4X_{t}-X_{t}^3)dt+X_{t}dW_{t},
\end{align}
where $W_{t}$ is a scalar Brownian motion. The generator of this system is known from (\ref{generator2}) by taking $b(x)=4x-x^3$ and $\sigma_{1}(x)=x$. \\
Thus we have
\begin{align}
(\mathcal{L_B}f)(x)=(4x-x^3)f^{'}(x)+\frac{1}{2}x^2f^{''}(x).
\end{align}
We use the EDMD algorithm to data, from (\ref{SDE2}), to learn their governing law. The data set has  $10^6$ initial points on $x\in  (-2,2)$ drawn from a uniform grid, which constitute $X$, and their positions after $\Delta t=0.01$, which constitute $Y.$ The sample paths were obtained by Euler-Maruyama method. The dictionary chosen is a polynomial with order up to 5. \\
Under this setting, we obtain the estimation of drift term and diffusion term. As can be seen from Table \ref{table:1dBM}, it is very close to the true parameters. Using the stochastic differential equation coefficients learned from the data, we obtained mean exit time and escape probability by finite differential method to solve the partial differential equations (\ref{eqn:MET3}) and (\ref{eqn:EP3}). Here we use the same grid to discretize the true system and the learned system. We can see from Fig.\ref{fig1dBMa} and Fig.\ref{fig1dBMb} that they are almost identical.\\

\begin{table}
\caption{Identified coefficients for 1-D double-well system driven by Brownian motion.}
\label{table:1dBM}
\centering
\subtable[Drift term]{
       \begin{tabular}{ccc}
        \hline
        basis& true& learning\\
        \hline
        $1$& 0& 0\\
        $x$& 4& 4.1243\\
        $x^2$& 0& 0\\
        $x^3$& -1& -1.1372\\
        $x^4$& 0& 0\\
        $x^5$& 0& 0\\
        \hline
       \end{tabular}
       \label{tab:BM1dDrift}
}
\qquad
\subtable[Diffusion term $\sigma_{1}^2$]{
       \begin{tabular}{ccc}
        \hline
        basis& true& learning\\
        \hline
        $1$& 0& 0\\
        $x$& 0& 0\\
        $x^2$& 1& 0.8823\\
        $x^3$& 0& 0\\
        $x^4$& 0& 0\\
        $x^5$& 0& 0\\
        \hline
       \end{tabular}
       \label{tab:BM1dDiffusion}
}
\end{table}

\begin{figure}
\label{FIG:1dBM}
\centering
\subfigure[Mean exit time]{
\label{fig1dBMa} 
\includegraphics[width=5cm]{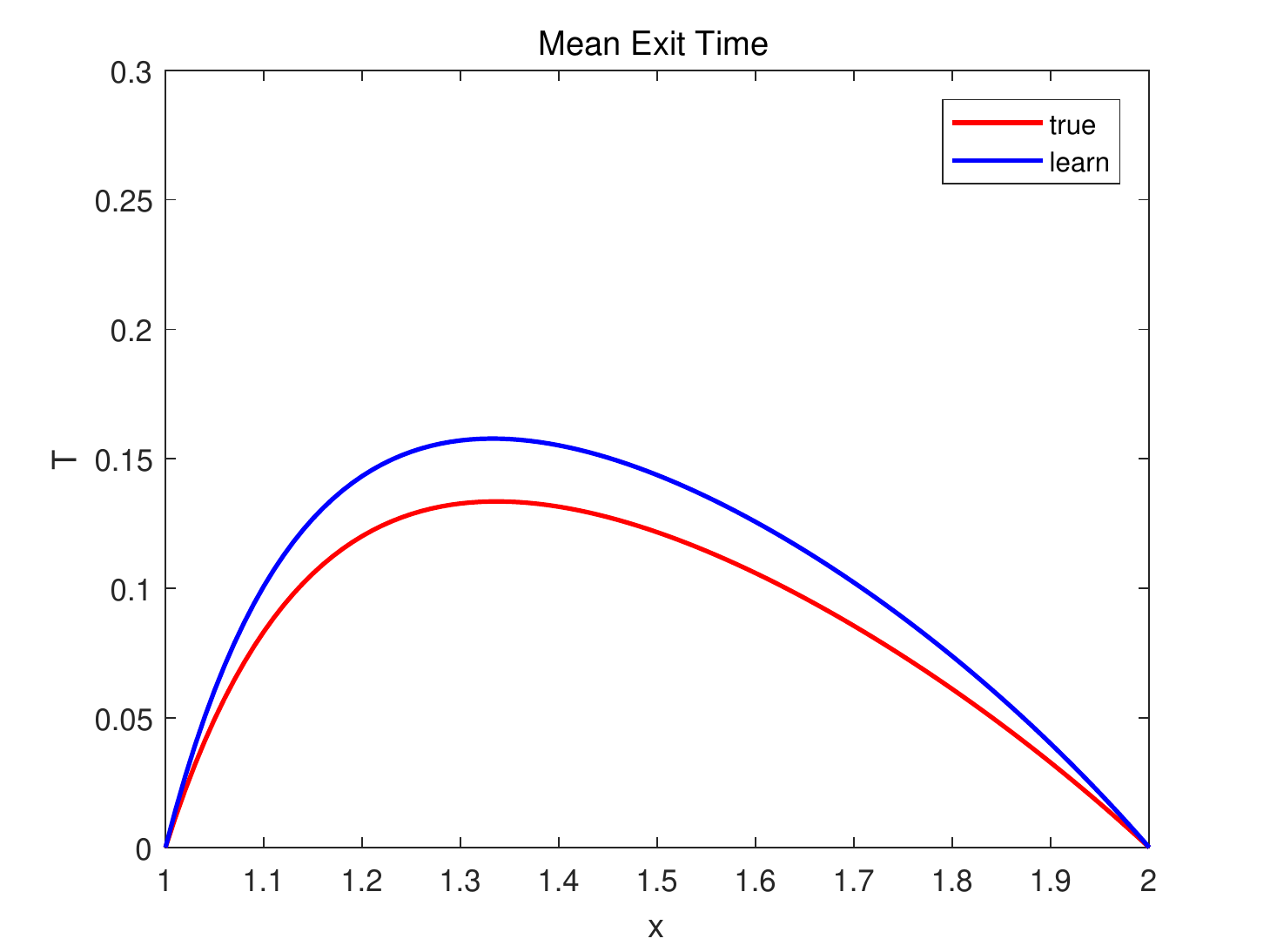}}
\hspace{0in}
\subfigure[Escape probability]{
\label{fig:subfig:b} 
\includegraphics[width=5cm]{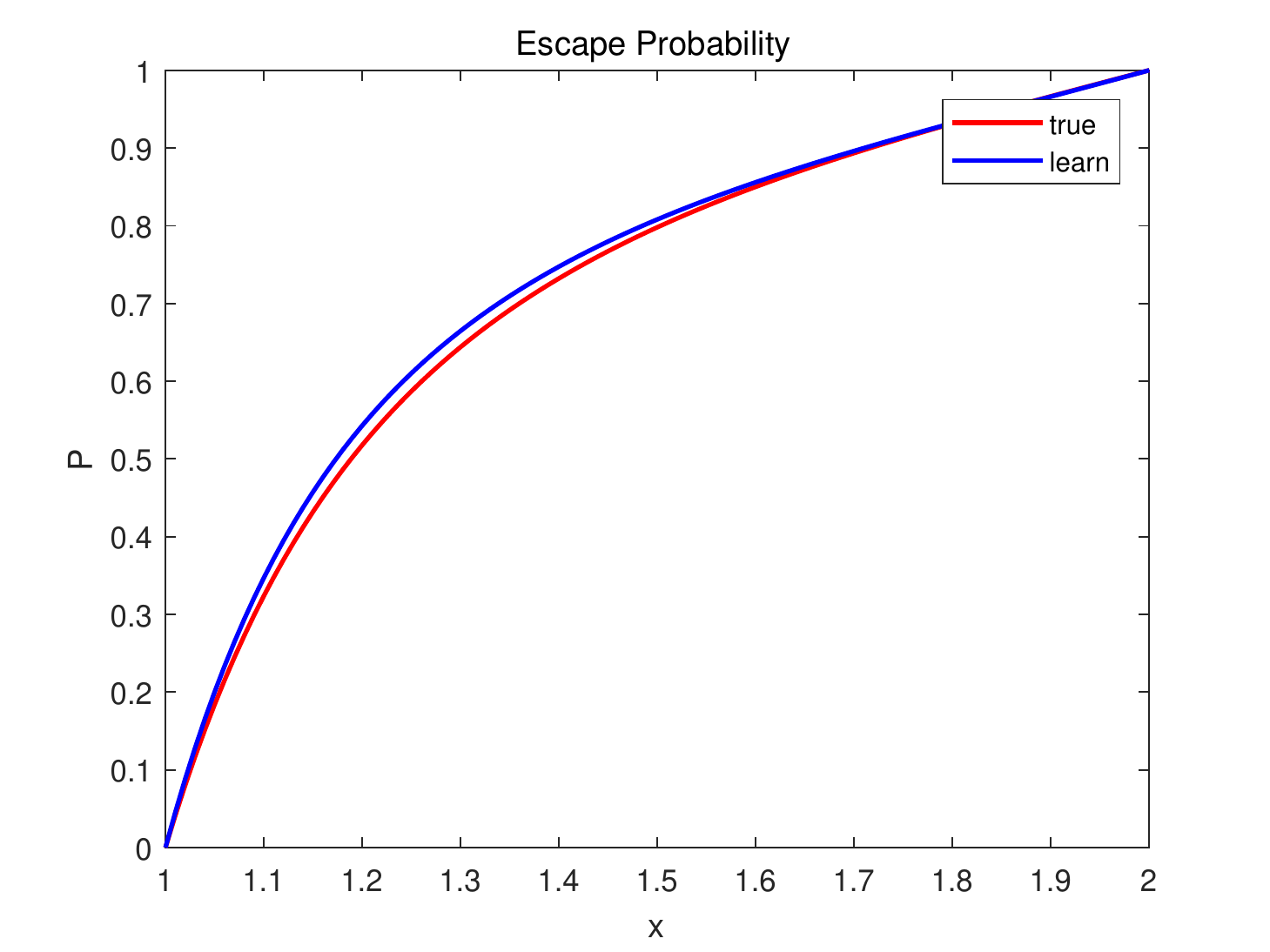}
\label{fig1dBMb}}
\caption{1-D double-well system driven by Brownian motion. (a) Mean exit time for the true system and the learned  system. (b) Escape probability for the true system and the learned system.}
\end{figure}
\end{Example}

\begin{Example} Consider a two-dimensional system with two-dimensional Brownian motion
\begin{align}\label{SDE3}
dX_{t}=(3X_{t}-Y_{t}^2)dt+X_{t}dW_{t}, \nonumber\\
dY_{t}=(2X_{t}+Y_{t})dt+Y_{t}d\tilde{W}_{t},
\end{align}
where $W_{t}$ and $\tilde{W}_{t}$ are two disjoint independent scalar real-value Brownian motion. The generator of this system  is easy to know from (\ref{generator2}) by taking $b_{1}(x,y)=3x-y^2$, $b_{2}(x,y)=2x+y$ and

\begin{center} {$\sigma_{1} = \left[ {\begin{array}{*{20}{c}}
x&0\\
0&y
\end{array}} \right].$}

\end{center}
Similar to Example 4.1, We use the EDMD algorithm to data, from (\ref{SDE3}), to learn their governing law. The data set has  $10^6$ initial points on $x\in  (-1,1)\times(-1,1)$ drawn from a uniform grid, which constitute $X$, and their positions after $\Delta t=0.01$, which constitute $Y.$ The sample paths were obtained by Euler-Maruyama method. The dictionary chosen is a polynomial with order up to 3. \\
As can be seen from Table \ref{table:2dBMdrift} and Table \ref{table:2dBMdiffusion}, it is very close to the true parameters. Using the stochastic differential equation coefficients learned from the data, we obtained mean exit time and escape probability by solving the partial differential equations (\ref{eqn:MET3}) and (\ref{eqn:EP3}). We can see from Fig.\ref{figMET2dBMa}, Fig.\ref{figMET2dBMb}, Fig.\ref{figEP2dBMa} and Fig.\ref{figEP2dBMb} that they are almost identical. The mean error of mean exit time and escape time are $0.0041$ and $0.0575$. These results show that we still can obtain the accurate mean exit time and escape probability in high dimension from data.\\

\begin{table}

\caption{Identified drift terms for 2-D system driven by Brownian motion.}
\label{table:2dBMdrift}
\centering
\subtable[Drift term $b_{1}$]{
       \begin{tabular}{ccc}
        \hline
        basis& true& learning\\
        \hline
        $1$& 0& 0\\
        $x$& 3& 2.9837\\
        $y$& 0& 0\\
        $x^2$& 0& 0\\
        $xy$& 0& 0\\
        $y^2$& -1&  -1.0026\\
        $x^3$& 0& 0\\
        $x^2y$& 0& 0\\
        $xy^2$& 0& 0\\
        $y^3$& 0& 0\\
        \hline
       \end{tabular}
       \label{tab:BM2dDrift}
}
\qquad
\subtable[Drift term $b_{2}$]{
       \begin{tabular}{ccc}
       \hline
        basis& true& learning\\
        \hline
        $1$& 0& 0\\
        $x$& 2& 2.0391\\
        $y$& 1& 1.0056\\
        $x^2$& 0& 0\\
        $xy$& 0& 0\\
        $y^2$& 0& 0\\
        $x^3$& 0& 0\\
        $x^2y$& 0& 0\\
        $xy^2$& 0& 0\\
        $y^3$& 0& 0\\
        \hline
       \end{tabular}
       \label{tab:BM2dDiffusion}
}
\end{table}

\begin{table}

\caption{Identified diffusion terms for 2-D system driven by Brownian motion.}
\label{table:2dBMdiffusion}
\centering
\subtable[Diffusion term $\sigma_{11}^2$]{
       \begin{tabular}{ccc}
        \hline
        basis& true& learning\\
        \hline
        $1$& 0& 0\\
        $x$& 0& 0\\
        $y$& 0& 0\\
        $x^2$& 1&  1.2835\\
        $xy$& 0& 0\\
        $y^2$& 0& 0\\
        $x^3$& 0& 0\\
        $x^2y$& 0& 0\\
        $xy^2$& 0& 0\\
        $y^3$& 0& 0\\
        \hline
       \end{tabular}
       \label{tab:BM2dDrift}
}
\qquad
\subtable[Diffusion term $\sigma_{22}^2$]{
       \begin{tabular}{ccc}
       \hline
        basis& true& learning\\
        \hline
        $1$& 0& 0\\
        $x$& 0& 0\\
        $y$& 0& 0\\
        $x^2$& 0& 0\\
        $xy$& 0& 0\\
        $y^2$& 1& 0.9750\\
        $x^3$& 0& 0\\
        $x^2y$& 0& 0\\
        $xy^2$& 0& 0\\
        $y^3$& 0& 0\\
        \hline
       \end{tabular}
       \label{tab:BM2dDiffusion}
}
\end{table}

\begin{figure}
\label{FIG:MET2dBM}
\centering
\subfigure[Mean exit time for the true system]{
\label{figMET2dBMa} 
\includegraphics[width=5cm]{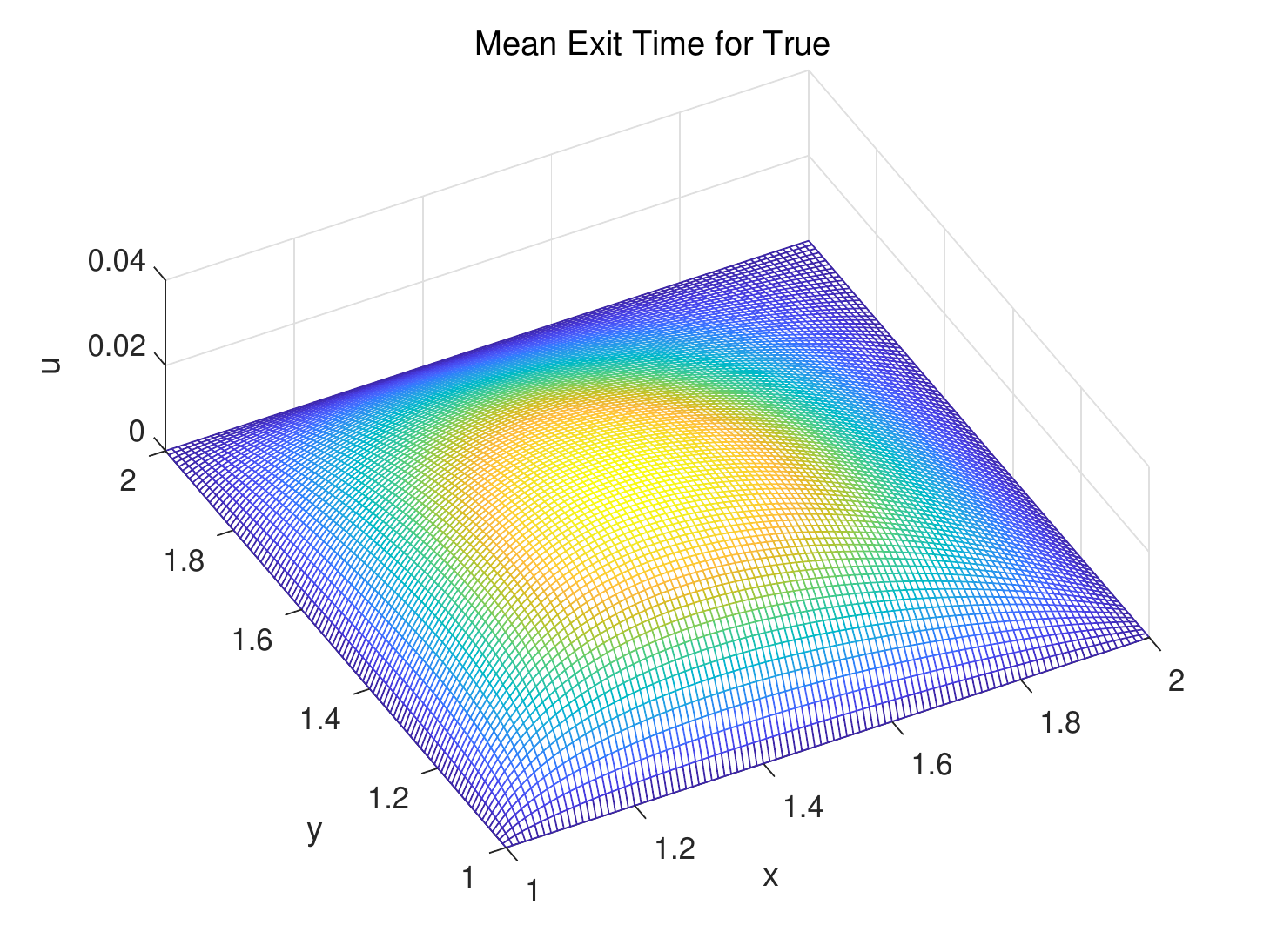}}
\hspace{0in}
\subfigure[Mean exit time for the learning system]{
\label{fig:subfig:b} 
\includegraphics[width=5cm]{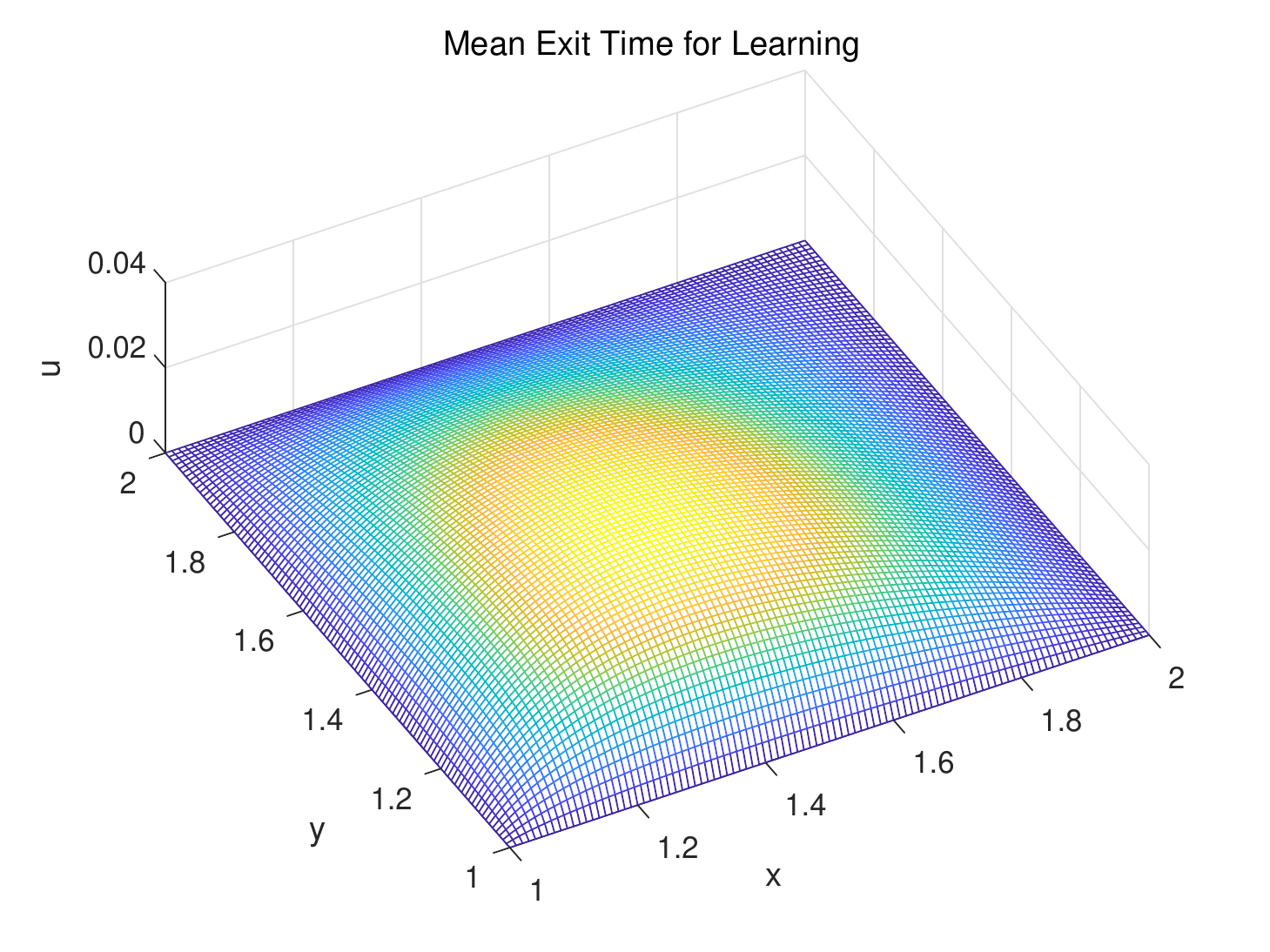}
\label{figMET2dBMb}}
\caption{2-D system driven by Brownian motion. (a) Mean exit time for the true system. (b) Mean exit time for the learned  system.}
\end{figure}

\begin{figure}
\label{FIG:EP2dBM}
\centering
\subfigure[Escape probability for the true system]{
\label{figEP2dBMa} 
\includegraphics[width=5cm]{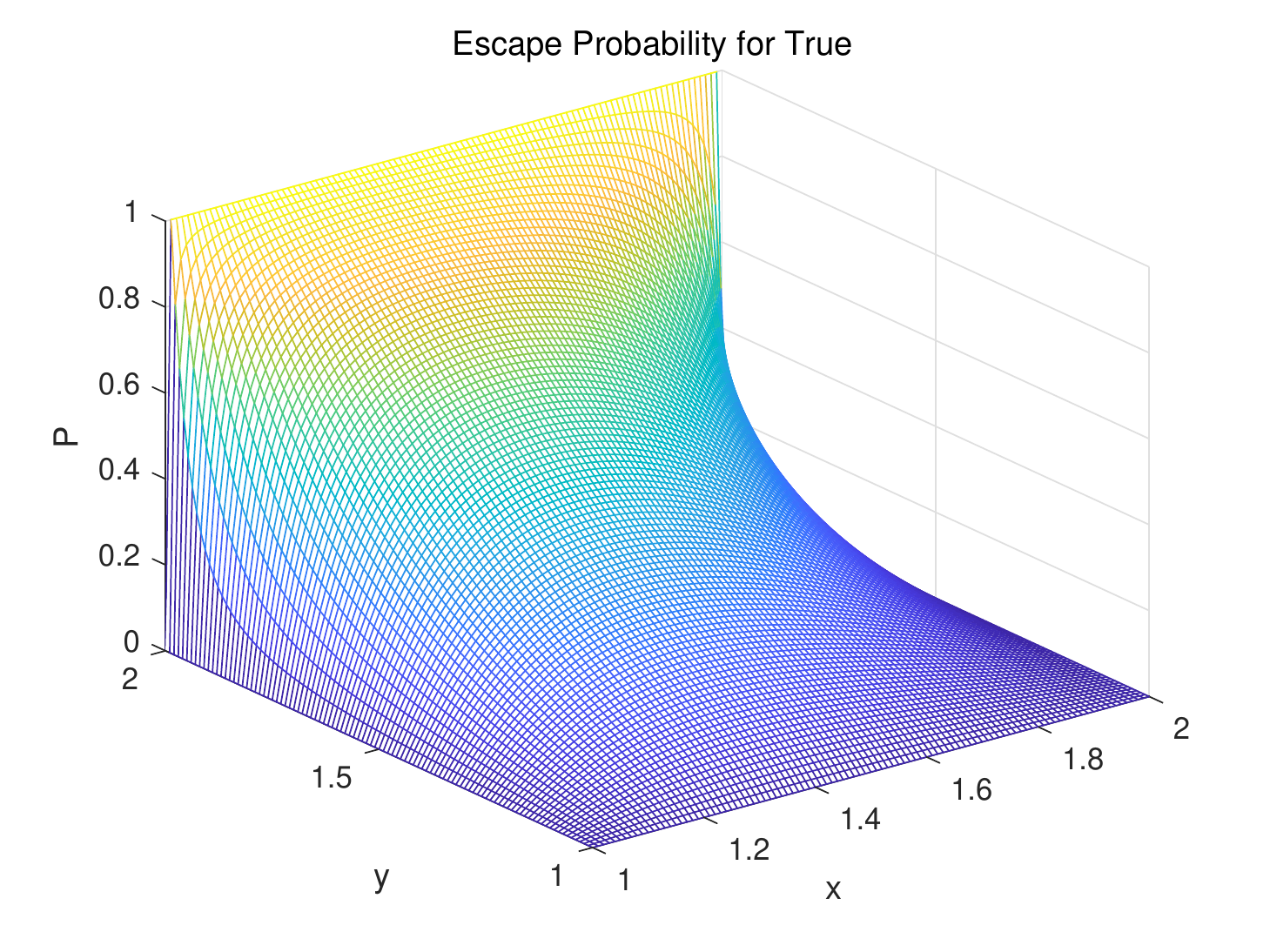}}
\hspace{0in}
\subfigure[Escape probability for the learning system]{
\label{fig:subfig:b} 
\includegraphics[width=5cm]{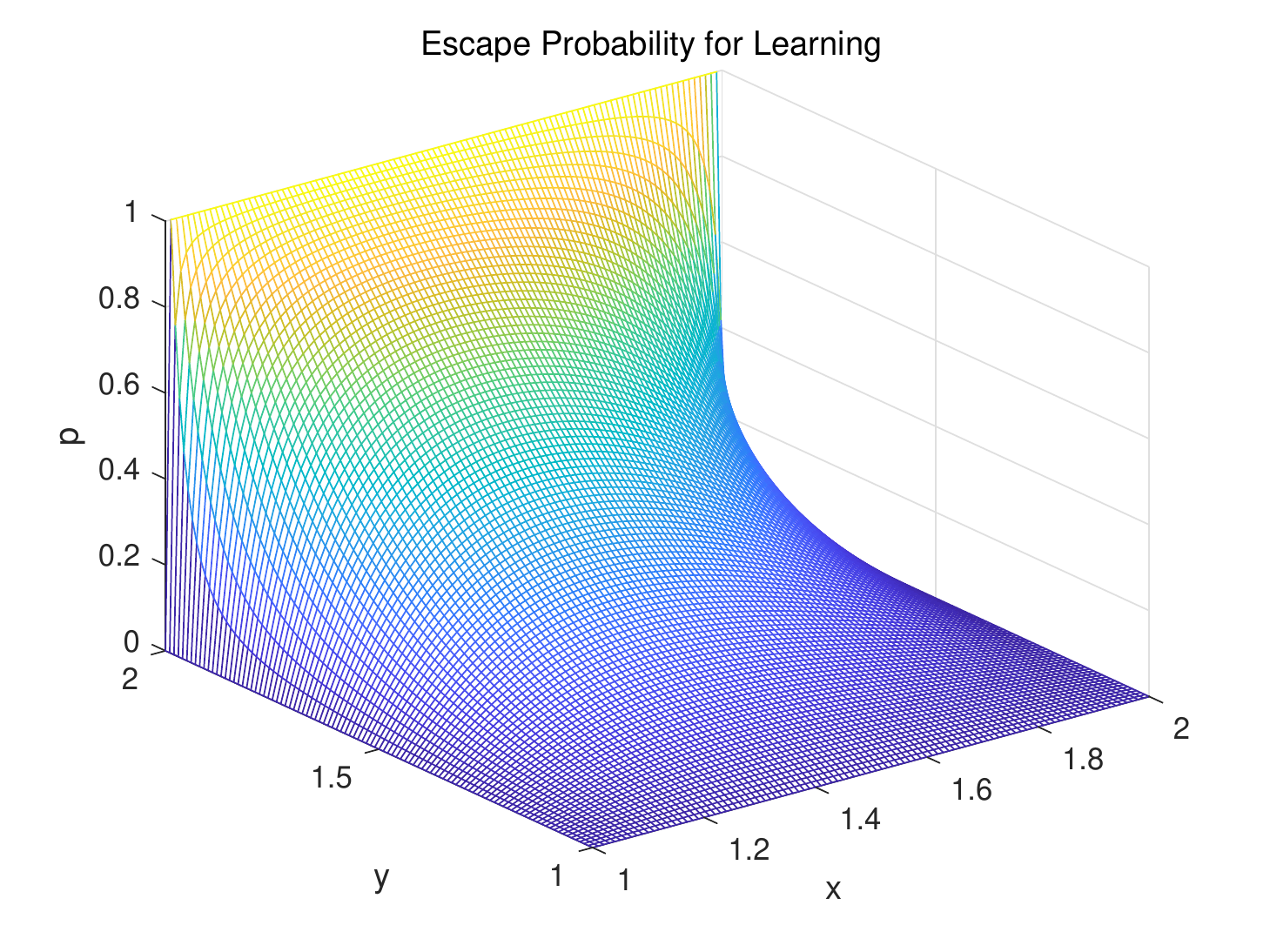}
\label{figEP2dBMb}}
\caption{2-D system driven by Brownian motion. (a) Escape probability for the true system. (b) Escape probability for the learned  system.}
\end{figure}
\end{Example}

\subsection{Dynamical systems with L\'evy motion  }

We will demonstrate that we can   obtain mean exit time  and escape probability from data, simulated   for systems with rotationally invariant $\alpha$-stable L\'evy motion. In what follows computations,  we  take that the stable index $\alpha=1$ and the jump measure $\tilde{\nu}(dx)=I_{\{\Vert x\Vert<c\}}\frac{dx}{|x|^{d+\alpha}}$ with $c=1$.\\
\\
\begin{Example} Consider a scalar double-well system with scalar L\'evy motion as follows
\begin{align}
dX_{t}=(4X_{t}-X_{t}^3)dt+X_{t}dW_{t}+d\tilde{L}_{t}^\alpha,
\end{align}
where $W_{t}$ is a scalar real-value Brownian motion and $L_{t}^\alpha$ is a L\'evy motion with triple $(0,0,\tilde{\nu})$ as we said in Section \ref{Koopman}. The generator of this system is easy to know from (\ref{generator1}) by taking $b(x)=4x-x^3$,$\sigma_{1}=x$ and $\sigma_{2}=1$. \\
Thus we have

\begin{align}
(\mathcal{L}f)(x)&=(4x-x^3)f^{'}(x)+\frac{1}{2}x^2f^{''}(x) \nonumber\\
 &+ C_{\alpha}\int_{\mathbb{R}\backslash\{0\}} [f(x+y)-f(x)]I_{\{\Vert y\Vert<c\}}(y)\, \frac{dy}{|y|^{1+\alpha}}.
\end{align}

The data set has  $10^6$ initial points on $x\in  (-2,2)$ drawn from a uniform grid, which constitute $X$, and their positions after $\Delta t=0.01$, which constitute $Y.$ The sample paths were obtained by Euler-Maruyama method. The dictionary chosen is a polynomial with order up to 5. \\
Notice that we have non-Gaussian noise here. As can be seen from Table \ref{table:1dLPdrift}, it is very close to the true parameters. In addition, the diffusion term squared of L\'evy motion, learning from data, is equal to 1.0955. It's close to the true coefficient 1. Using the stochastic differential equation coefficients learned from the data, we obtained mean exit time and escape probability by solving the partial differential equations (\ref{eqn:MET1}), (\ref{eqn:MET2}) and (\ref{eqn:EP}). We can see from Fig.\ref{fig1dLPa} and Fig.\ref{fig1dLPb} that they are almost identical.\\

\begin{table}

\caption{Identified coefficients for 1-D double-well system driven by L\'evy motion.}
\label{table:1dLPdrift}
\centering
\subtable[Drift term]{
       \begin{tabular}{ccc}
        \hline
        basis& true& learning\\
        \hline
        $1$& 0& 0\\
        $x$& 4& 4.0716\\
        $x^2$& 0& 0\\
        $x^3$& -1& -1.0966\\
        $x^4$& 0& 0\\
        $x^5$& 0& 0\\
        \hline
       \end{tabular}
       \label{tab:BM1dDrift}
}
\qquad
\subtable[Diffusion term $\sigma_{1}^2$]{
       \begin{tabular}{ccc}
        \hline
        basis& true& learning\\
        \hline
        $1$& 0& 0\\
        $x$& 0& 0\\
        $x^2$& 1& 0.9862\\
        $x^3$& 0& 0\\
        $x^4$& 0& 0\\
        $x^5$& 0& 0\\
        \hline
       \end{tabular}
       \label{tab:BM1dDiffusion}
}
\end{table}

\begin{figure}
\label{FIG:1dLP}
\centering
\subfigure[Mean exit time]{
\label{fig1dLPa} 
\includegraphics[width=5cm]{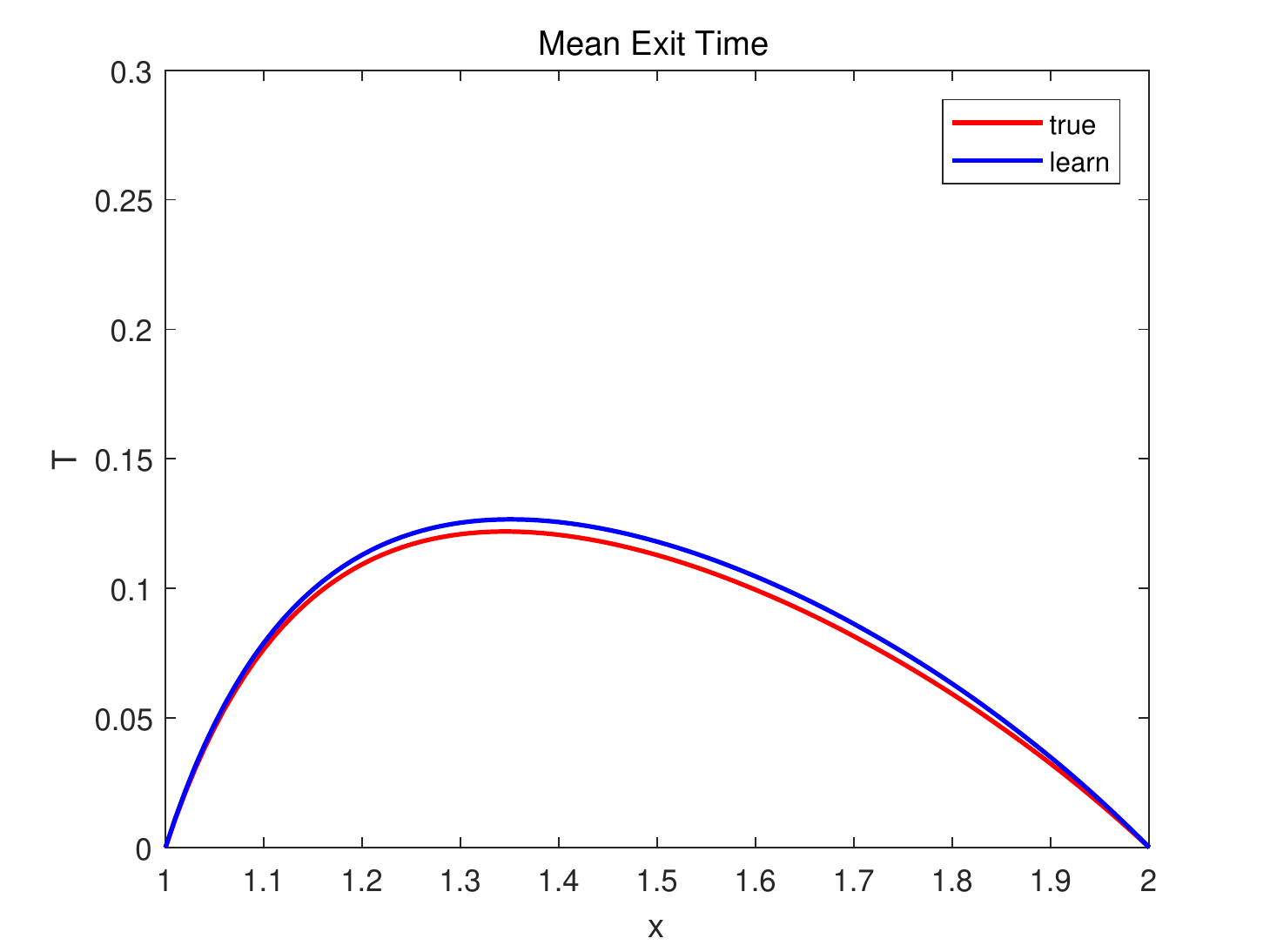}}
\hspace{0in}
\subfigure[Escape probability]{
\label{fig:subfig:b} 
\includegraphics[width=5cm]{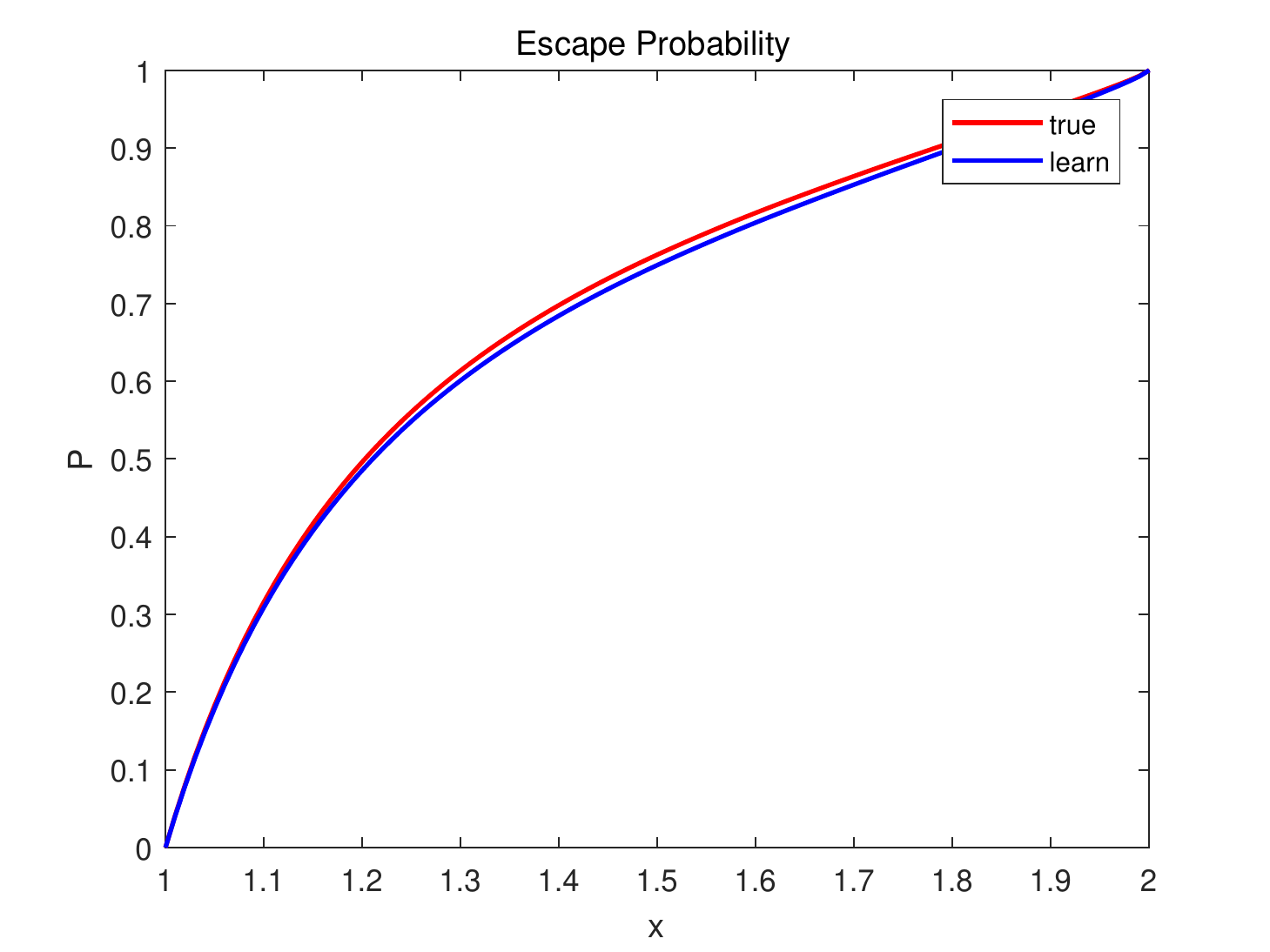}
\label{fig1dLPb}}
\caption{1-D double-well system driven by L\'evy motion. (a) Mean exit time for the true system and the learned system. (b) Escape probability for the true system and the learned  system.}
\end{figure}
\end{Example}

\begin{Example} Consider a two-dimensional system
\begin{align}
dX_{t}=(3X_{t}-Y_{t}^2)dt+X_{t}dW_{t}+d\tilde{L}_{1,t}^\alpha, \nonumber\\
dY_{t}=(2X_{t}+Y_{t})dt+Y_{t}d\tilde{W}_{t}+d\tilde{L}_{2,  t}^\alpha,
\end{align}
where $W_{t}$, $\tilde{W}_{t}$, $\tilde{L}_{1,t}$ and $\tilde{L}_{2,t}$ are four disjoint independent scalar real-value Brownian motions and L\'evy motions as we said in Section \ref{Koopman}. Note that we still assume that the jump size of L\'evy motions are bounded. The generator of this system is easy to know from (\ref{generator1}) by taking $b_{1}(x,y)=3x-y^2$, $b_{2}(x,y)=2x+y$ and

\begin{center} {$\sigma_{1} = \left[ {\begin{array}{*{20}{c}}
x&0\\
0&y
\end{array}} \right],$ $\sigma_{2} = \left[ {\begin{array}{*{20}{c}}
1&0\\
0&1
\end{array}} \right].$}
\end{center}

For this 2-D system, the data set has  $10^6$ initial points on $x\in  (-1,1)\times(-1,1)$ drawn from a uniform grid, which constitute $X$, and their positions after $\Delta t=0.01$, which constitute $Y.$ The sample paths were obtained by Euler-Maruyama method. The dictionary chosen is a polynomial with order up to 3. \\ As can be seen from Table \ref{table:2dLPdrift} and Table \ref{table:2dLPdiffusion}, it is very close to the true parameters. In addition, the diffusion term squared of L\'evy motion, learning from data, is equal to $[1.0710, 0; 0, 1.0370]$. It's close to the true coefficient $[1,0; 0, 1]$. Using the stochastic differential equation coefficients learned from the data, we obtained mean exit time and escape probability by solving the partial differential equations (\ref{eqn:MET1}), (\ref{eqn:MET2}) and (\ref{eqn:EP}). We can see from Fig.\ref{figMET2dLPa}, Fig.\ref{figMET2dLPb}, Fig.\ref{figEP2dLPa} and Fig.\ref{figEP2dLPb} that they are almost identical. The mean error of mean exit time and escape time are $6.2429e-04$ and $0.0128$.\\

\begin{Remark}
These examples show that although we use finite dimensional matrices to approximate infinite dimensional linear operators, we can still accurately capture the mean exit time and escape probability of the corresponding stochastic dynamical system. i.e., we can discover deterministic indexes (e.g., mean exit time, escape probability) to gain insights about transition phenomena from data for stochastic dynamical systems under Gaussian noise or non-Gaussian noise.
\end{Remark}
\begin{table}

\caption{Identified drift terms for 2-D system driven by L\'evy motion.}
\label{table:2dLPdrift}
\centering
\subtable[Drift term $b_{1}$]{
       \begin{tabular}{ccc}
        \hline
        basis& true& learning\\
        \hline
        $1$& 0& 0\\
        $x$& 3& 3.1096\\
        $y$& 0& 0\\
        $x^2$& 0& 0\\
        $xy$& 0& 0\\
        $y^2$& -1& -0.9587\\
        $x^3$& 0& 0\\
        $x^2y$& 0& 0\\
        $xy^2$& 0& 0\\
        $y^3$& 0& 0\\
        \hline
       \end{tabular}
       \label{tab:BM2dDrift}
}
\qquad
\subtable[Drift term $b_{2}$]{
       \begin{tabular}{ccc}
       \hline
        basis& true& learning\\
        \hline
        $1$& 0& 0\\
        $x$& 2& 2.0882\\
        $y$& 1& 1.0324\\
        $x^2$& 0& 0\\
        $xy$& 0& 0\\
        $y^2$& 0& 0\\
        $x^3$& 0& 0\\
        $x^2y$& 0& 0\\
        $xy^2$& 0& 0\\
        $y^3$& 0& 0\\
        \hline
       \end{tabular}
       \label{tab:BM2dDiffusion}
}
\end{table}

\begin{table}

\caption{Identified drift terms for 2-D system driven by L\'evy motion.}
\label{table:2dLPdiffusion}
\centering
\subtable[Diffusion term $\sigma_{11}^2$]{
       \begin{tabular}{ccc}
        \hline
        basis& true& learning\\
        \hline
        $1$& 0& 0\\
        $x$& 0& 0\\
        $y$& 0& 0\\
        $x^2$& 1& 0.9668\\
        $xy$& 0& 0\\
        $y^2$& 0& 0\\
        $x^3$& 0& 0\\
        $x^2y$& 0& 0\\
        $xy^2$& 0& 0\\
        $y^3$& 0& 0\\
        \hline
       \end{tabular}
       \label{tab:BM2dDrift}
}
\qquad
\subtable[Diffusion term $\sigma_{22}^2$]{
       \begin{tabular}{ccc}
       \hline
        basis& true& learning\\
        \hline
        $1$& 0& 0\\
        $x$& 0& 0\\
        $y$& 0& 0\\
        $x^2$& 0& 0\\
        $xy$& 0& 0\\
        $y^2$& 1& 1.0213\\
        $x^3$& 0& 0\\
        $x^2y$& 0& 0\\
        $xy^2$& 0& 0\\
        $y^3$& 0& 0\\
        \hline
       \end{tabular}
       \label{tab:BM2dDiffusion}
}
\end{table}

\begin{figure}

\label{FIG:MET2dLP}
\centering
\subfigure[Mean exit time for the true system]{
\label{figMET2dLPa} 
\includegraphics[width=5cm]{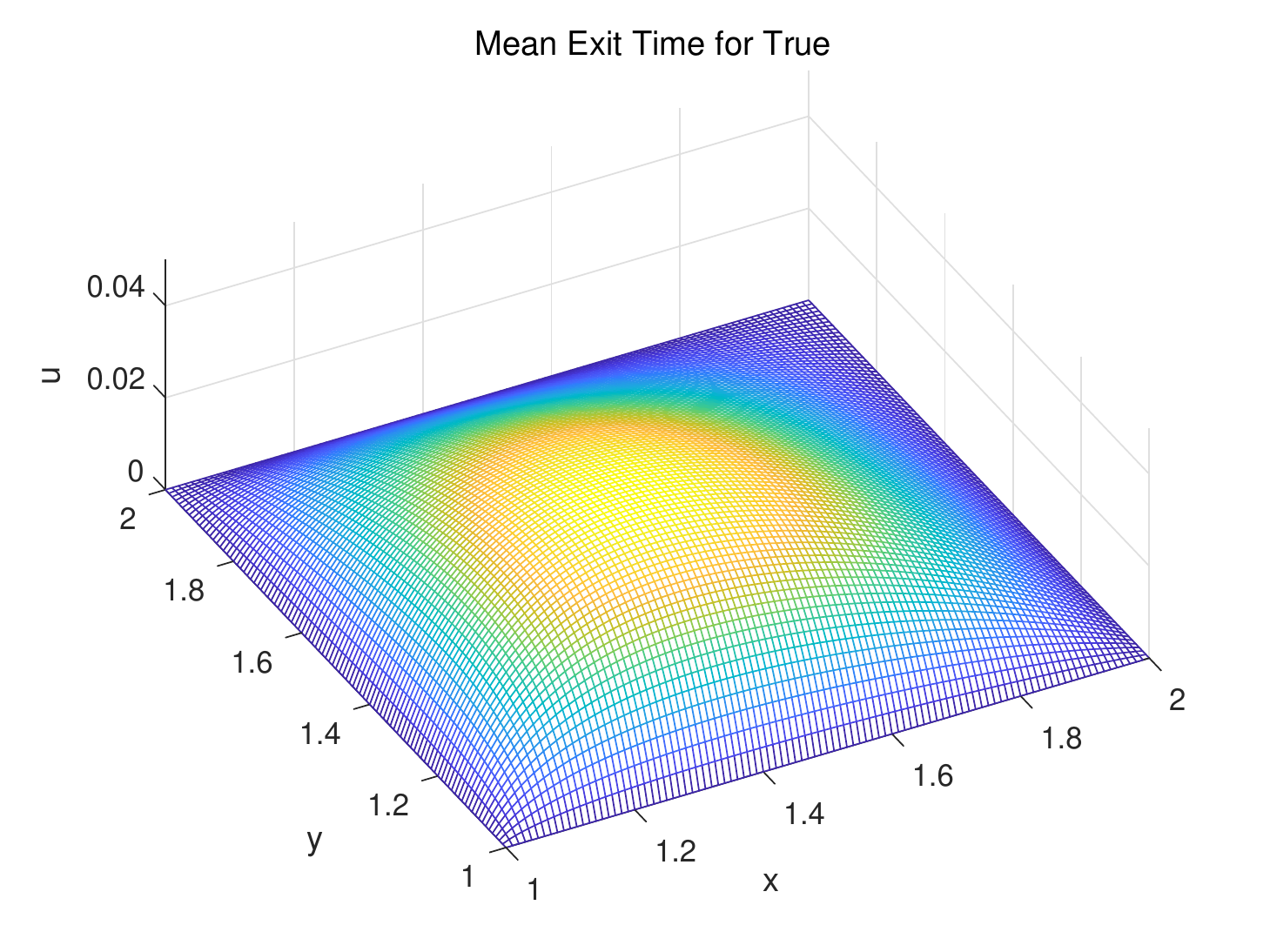}}
\hspace{0in}
\subfigure[Mean exit time for the learning system]{
\label{fig:subfig:b} 
\includegraphics[width=5cm]{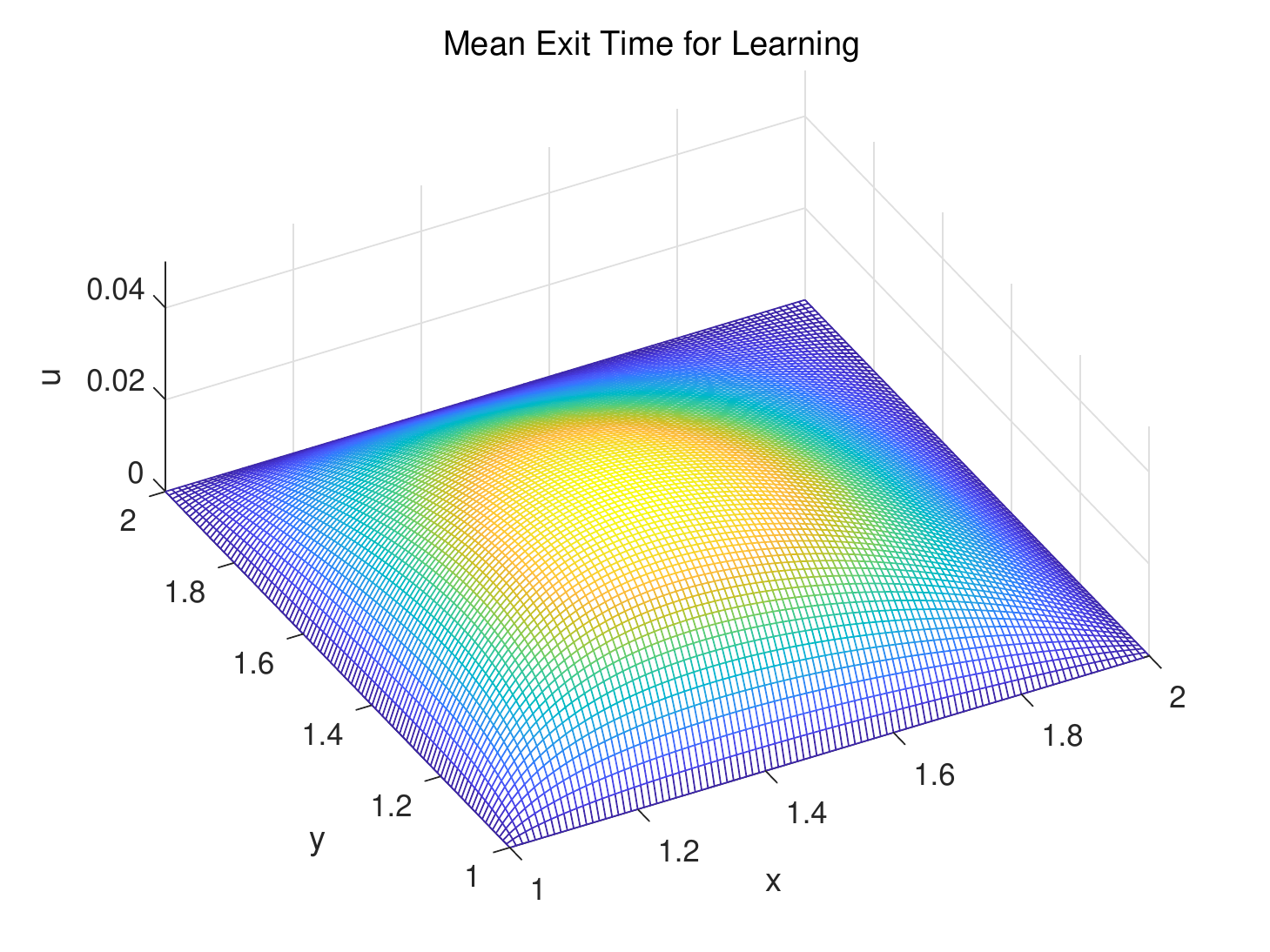}
\label{figMET2dLPb}}
\caption{2-D system driven by L\'evy motion. (a) Mean exit time for the true system. (b) Escape probability for the learned  system.}
\end{figure}

\begin{figure}
\label{FIG:EP2dLP}
\centering
\subfigure[Escape probability for the true system]{
\label{figEP2dLPa} 
\includegraphics[width=5cm]{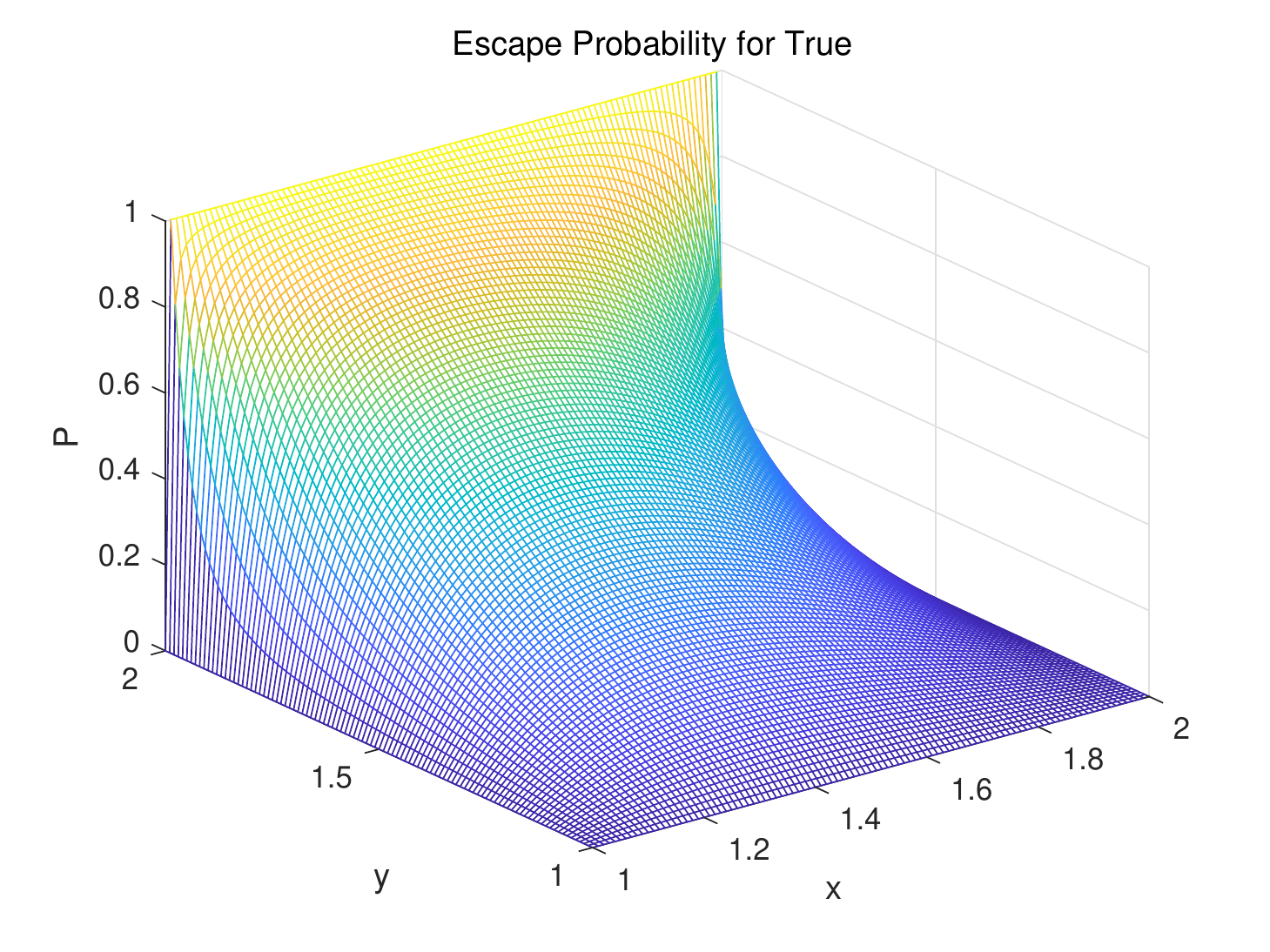}}
\hspace{0in}
\subfigure[Escape probability for the learning system]{
\label{fig:subfig:b} 
\includegraphics[width=5cm]{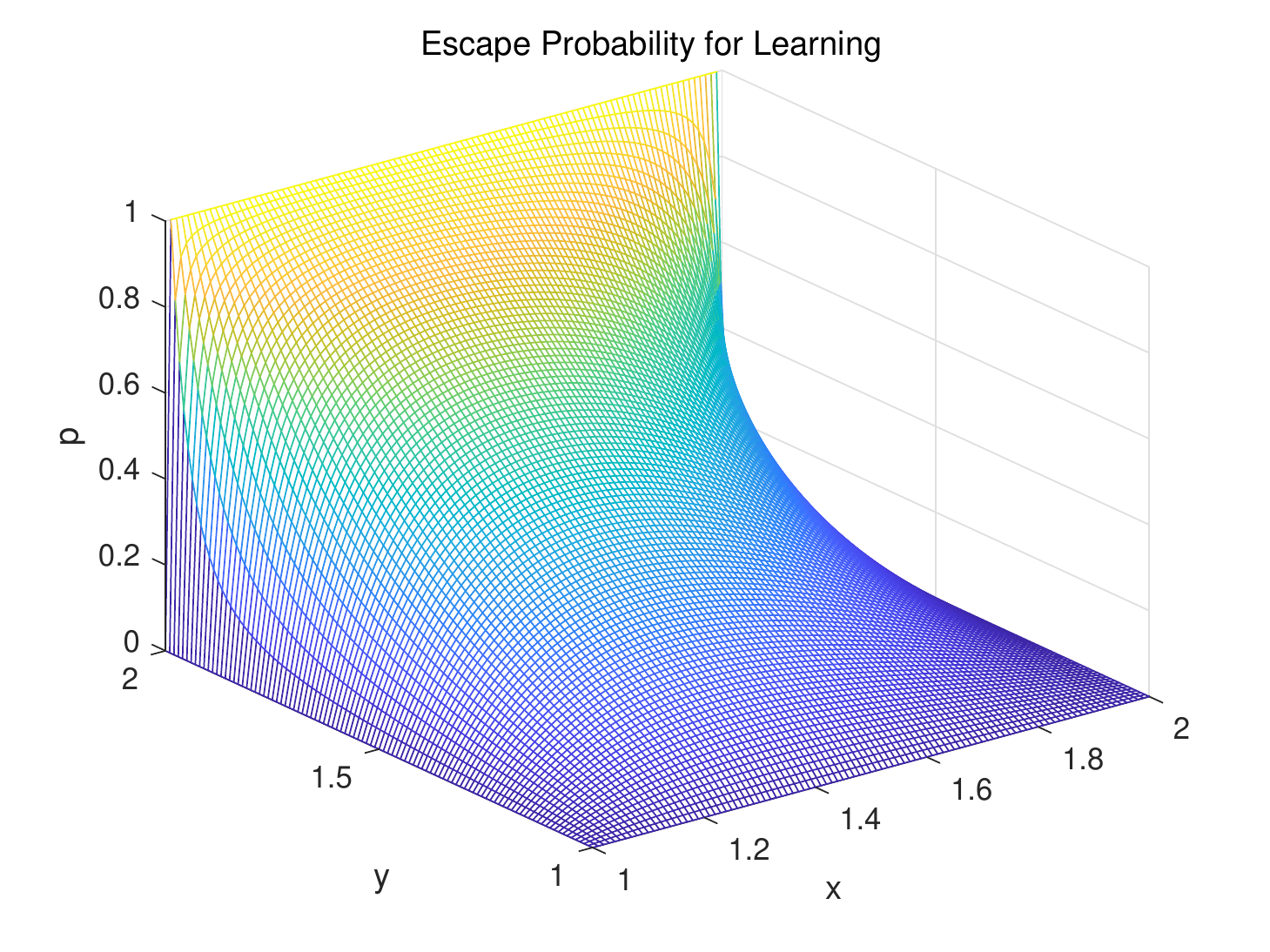}
\label{figEP2dLPb}}
\caption{2-D system driven by L\'evy motion. (a) Escape probability for the true system. (b) Escape probability for the learned  system.}
\end{figure}
\end{Example}

\section{Discussion}

We have presented   an approach  for obtaining deterministic indexes about  transition phenomena of  non-Gaussian stochastic dynamical systems from data.  In particular,  we have identified the coefficients from the stochastic differential equations with non-Gaussian L\'evy motion and extracted information about   transition phenomena (e.g. mean exit time and escape probability) from the data of the stochastic system
.   Here we have applied the Koopman analysis framework to deal with stochastic systems with L\'evy motion.

Existing relevant works  focused on either deterministic systems,  or stochastic systems with Gaussian Brownian motion, including  our earlier work \cite{WuFuDuan}, where the stochastic governing law is directly learned from data.


In fact, we   can further obtain transition probability density functions of the stochastic dynamical systems, because we can additionally solve the associated Fokker-Planck equations.

However, there are several challenges for further investigation. First, so far there is  no rigorous error analysis for   the finite-dimensional approximation of infinite dimensional linear operators. Generally speaking, finite-dimensional approximation is more effective when operators have pure point spectra. However, when the linear operator has a continuous spectrum, there is still no good way to deal with it. Second,  how to choose  an appropriate  basis of functions is a challenge. Third,   we have noticed  that stochastic dynamical systems need more simulated data than deterministic dynamical systems, in order to extract transition dynamics.  Hence,  reducing the demand for   data  is   worth of further investigation.  Fourth,  our method requires data on system states. If we only have observation and this observation function can be expanded using basis, our method still can work. Otherwise, we have to use  other suitable methods based on   observations. For example, our earlier work  \cite{Gao1} demonstrated a method to learn stochastic governing laws with observations not on system states. Finally, our method is not valid when the L\'evy noise is   non-diagonal matrix or multiplicative, because we can not separate non-Gaussian noise from Gaussian noise.

\section*{Acknowledgements}
The authors would like to thank Ting Gao,   Min Dai,   Qiao Huang,   Xiaoli Chen,   Pingyuan Wei,  and  Yancai Liu for helpful discussions. This work was partly supported by the   NSFC grants 11531006 and 11771449.

\section*{Data Availability Statements}
The data that support the findings of this study are openly available in GitHub, reference number \cite{Code}.

\section*{Appendix}\label{secAppendix}
\textbf{L\'evy motions} \\
Let $L=(L_{t},t\geq0)$ be a stochastic process defined on a probability space $(\Omega, \mathcal{F}, P)$. We say that L is a L\'evy motion if:
\begin{itemize}
\item $L_{0}=0 (a.s.)$;
\item $L$ has independent and stationary increments;
\item $L$ is stochastically continuous, i.e., for all $a>0$ and for all $s\geq0$
$$
\lim\limits_{t\to s}P(|L_{t}-L_{s}|>a)=0.
$$
\end{itemize}

\textbf{Characteristic functions}\\
For a L\'evy motion $(L_{t},t\geq0)$, we have the L\'evy-Khinchine formula,
$$
\mathbb{E}[e^{i(u,L_t)}] = exp\{t[i(b,u)-\frac{1}{2}(u,Au)+\int_{\mathbb{R}^d\backslash\{0\}} [e^{i(u,y)}-1-i(u,y)I_{\{\Vert y\Vert<1\}}(y)]\, \nu(dy)]\},
$$
for each $t\geq0$, $u\in\mathbb{R}^d$, where $(b,A,\nu)$ is the triple of L\'evy motion $(L_{t}, t\geq0)$.\\

\textbf{Theorem (The L\'evy-It\^{o} decomposition)} If $(L_{t},t\geq0)$ is a L\'evy motion with $(b,A,\nu)$, then there exists $b\in\mathbb{R}^d$, a Brownian motion $B_{A}$ with covariance matrix $A$ and an independent Poisson random measure $N$ on $\mathbb{R}^{+}\times(\mathbb{R}^d-\{0\})$ such that, for each $t\geq0$, \\
$$
L_{t}=bt+B_{A}(t)+\int_{0<|x|<c} x\, \tilde{N}(t,dx)+\int_{|x|\geq c} x\, N(t,dx),
$$
where $\int_{|x|\geq c} x\, N(t,dx)$ is a Poisson integral and $\int_{0<|x|<c} x\, \tilde{N}(t,dx)$ is a compensated Poisson integral defined by
$$
\int_{0<|x|<c} x\, \tilde{N}(t,dx)=\int_{0<|x|<c} x\, N(t,dx)-t\int_{0<|x|<c} x\, \nu(dx).
$$
Specially, in our work, we assume that the components of a $d-$dimensional L\'evy motion have triple $(0,0,\tilde{\nu})$. i.e., $$(\tilde{L}_{t}^{\alpha})_{i}=\int_{0<|x|<c} x\, \tilde{N}(t,dx),$$
where $i=1, 2, \ldots, d$, $x\in \mathbb{R}$ and $\tilde{\nu}(dx)=C_{\alpha}I_{\{\Vert x\Vert<c\}}(x)\frac{dx}{|x|^{1+\alpha}}$.
See \cite{Applebaum,Duan}.\\

\end{document}